\def\softd{{\leavevmode\setbox1=\hbox{d}%
\hbox to 1.05\wd1{d\kern-0.4ex{\char039}\hss}}}
\def\softt{{\leavevmode\setbox1=\hbox{t}%
\hbox to \wd1{t\kern-0.6ex{\char039}\hss}}}
\def\softl{l\kern-0.45ex\raise0.1ex\hbox{'}\kern-0.10ex}
\def\softL{L\kern-0.8ex\raise0.1ex\hbox{'}\kern0.1ex}

\documentstyle[a4]{article}
\begin{document}
\vspace*{40mm}
\begin{center}
        {\Large {\bf{THE MODELLING AND ANALYSIS OF FRACTIONAL-ORDER
                     CONTROL SYSTEMS IN DISCRETE DOMAIN}}} \\
         \vspace{10mm}
         \normalsize
         {Ivo PETR\'A\v{S}, \softL{}ubom\'{\i}r DOR\v{C}\'AK, Imrich 	
	   KO\v{S}TIAL \\
               \vspace{1mm}
                {   Department of Informatics and Process Control  \\
                    BERG Faculty, Technical University of Ko\v{s}ice   \\
                    B. N\v{e}mcovej 3, 042 00 Ko\v{s}ice, Slovak Republic  
		    \\
                    phone:      (+42195) 6025172                   \\
                    e-mail: {\it\{petras, dorcak, kostial\}@tuke.sk}}   \\
                  }

\end{center}
\vskip 5mm
\begin{abstract}
\hspace{-6.5mm}
This paper deals with fractional-order controlled systems and fractional-order controllers in the discrete domain. The mathematical description by the fractional difference equations and properties of these systems are presented. A practical example for modelling the fractional-order control loop is shown and obtained results are discussed in conclusion.

\end{abstract}

\hspace{2.5mm}

\section*{1. INTRODUCTION}

Fractional calculus was used for modelling  of physical systems,
but we can find only few works dealing with the application of
this mathematical tool in control theory (e.g. \cite{Axtell, Dorcak,
Matignon, Petras1, Oustaloup, Podlubny}). These works used continuous
mathematical models of fractional order. The fractional-order
systems have a unlimited memory, being integer-order systems
cases in which the memory is limited.

It is necessary to realise that the case of discrete
fractional-order systems is very important for their description
to have a finite difference equation. Such an equation can be
obtained by numerical approximation and by the inverse Z-transform
of a~discrete transfer function.

The aim of this paper is to show, how by using the fractional
calculus, we can obtain a more general fractional-order models
of the controlled objects and general structure for the classical
$PID$ controller in the discrete domain.

\section*{2. DISCRETE FRACTIONAL CALCULUS}

The idea of fractional calculus has been known since the development
of the regular calculus, with the first reference
probably being associated with Leibniz and L'Hospital in 1695.

Fractional calculus is a generalisation of integration and differentiation
to non-integer order fundamental operator $_{a}D^{\alpha}_{t}$,
where $a$ and $t$ are the limits of the operation.
The continuous integro-differential operator is defined as
$$
 _aD^{\alpha}_{t} = \left \{
        \begin{array}{ll}
                \frac{d^{\alpha}}{dt^{\alpha}} & \mbox{$\Re(\alpha)>0$,} \\
                 1 & \mbox{$\Re(\alpha)=0$,} \\
                \int_{a}^{t} (d\tau)^{-\alpha} & \mbox{$\Re(\alpha)<0$.}
        \end{array}
        \right.
$$
The two definitions used for the general fractional differintegral
are the Gr\"unwald-Letnikov (GL) definition and the Riemann-Liouville (RL)
definition \cite{Oldham}. The GL is given here
\begin{equation}\label{GLD_d}
     _{a}D^{\alpha}_{t}f(t)=\lim_{h \to 0} h^{-\alpha}
     \sum_{j=0}^{[\frac{t-a}{h}]}(-1)^j {\alpha \choose j}f(t-jh),
\end{equation}
where $[x]$ means the integer part of $x$. The RL definition is given as
\begin{equation}\label{LRL}
   _{a}D_{t}^{\alpha}f(t)=
    \frac{1}{\Gamma (n -\alpha)}
    \frac{d^{n}}{dt^{n}}
    \int_{a}^{t}
    \frac{f(\tau)}{(t-\tau)^{\alpha - n + 1}}d\tau, \\
\end{equation}
for $(n-1 < \alpha <n)$ and
where $\Gamma (.)$ is the well known Euler's {\it Gamma}  function.

In general, the approximation for fractional operator of
order $\alpha$ can be expressed by generating function
$\omega(\xi^{-1})$, where $\xi^{-1}$ is the shift operator. This
generating function and its expansion determine both the form of
the approximation and the coefficients \cite{Lubich}.

For discrete approximation of the time derivative,  we can use
the generating function corresponding to the Z-transform of
backward difference rule,
\begin{equation}
   \omega(z^{-1}) = \frac{1-z^{-1}}{T},
\end{equation}
and performing the power series expansion (PSE) of $(1-z^{-1})^{\pm \alpha}$,
we obtain the Z version of the GL formula by using the short
memory principle \cite{Podlubny}, for the discrete equivalent
of the fractional-order integro-differential operator
$\omega(z^{-1})^{\pm \alpha}$,
\begin{equation}\label{DFD}
(\omega(z^{-1}))^{\pm \alpha} = T^{\mp \alpha} \sum_{j=0}^{[\frac{L}{T}]} (-1)^j
{\pm \alpha \choose j} z^{[\frac{L}{T}] - j},
\end{equation}
where $T$ is the sample period, $L$ is the memory length, $[x]$ is the
integer part of $x$ and $(-1)^j {\pm \alpha \choose j}$ are 
a~binomial coefficients $c_j^{(\alpha)}, \,(j=0, 1, \dots)$.
For its calculation we can use the following expression:
\begin{equation}\label{b_k}
  c_0^{(\alpha)} = 1, \qquad c_j^{(\alpha)} = \left (1 -
                       \frac{1+ (\pm \alpha)}{j}\right)c_{j-1}^{(\alpha)}.
\end{equation}

Another possibility for the discrete approximation is the use of
the trapezoidal (Tustin) rule as a generating function for PSE.  For the discrete equivalent of the fractional-order integro-differential operator,
we can write a general formula \cite{Matignon}
\begin{equation}
(\omega(z^{-1}))^{\pm \alpha} = \left ( \frac{2}{T}
\frac{1-z^{-1}}{1+z^{-1}} \right )^{\pm \alpha}
\end{equation}
for obtaining the coefficients and the form of the approximation.

A detailed review of the approximation methods (Calson's,
Chareff's, CFE, Matsuda's, \dots) for continuous and
discrete fractional-order models was done in work \cite{Vinagre}.

\section*{3. FRACTIONAL-ORDER CONTROL CIRCUIT}

We will be studying the control system shown in Fig.1,
where $G_{c}(z)$ is the controller transfer function,
$G_{s}(z)$ is the controlled system transfer function,
$W(z)$ is an input, $E(z)$ is an error, $U(z)$ is the output
from controller and $Y(z)$ is the output from system \cite{Dorf}.

   \vspace*{4.1cm}
   \vspace{-20mm}
   \begin{picture}(110,42)
\setlength{\unitlength}{1mm}
   \put(26,27){$W(z) \ \ +$}
   \put(37,19){$ -$}
   \put(49.5,27){$E(z)$}
   \put(81.0,27){$U(z)$}
   \put(117.0,27){$Y(z)$}
   \put(26,25){\vector(1,0){15.9}}
   \put(45,10){\vector(0,1){11.7}}
   \put(45,25){\circle{6}}
   \put(42.9,27.0){\line(1,-1){4.0}}
   \put(42.9,23.0){\line(1, 1){4.0}}
   \put(48,25){\vector(1,0){12.2}}
   \put(45,10){\line(1,0){72.1}}
   \put(60,20){\framebox(20,10) [cc]{$G_c(z)$}}
   \put(80,25){\vector(1,0){10.0}}
   \put(90,20){\framebox(20,10) [cc]{$G_s(z)$}}
   \put(110,25){\vector(1,0){15.9}}
   \put(117,10){\line(0,1){15.1}}
\end{picture}
   \vspace{-10mm}
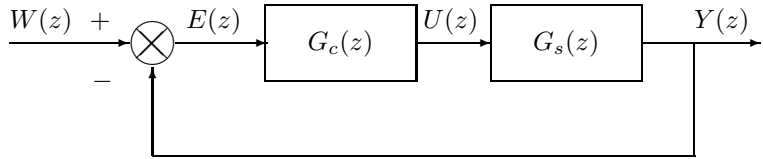
   \centerline{Figure 1: {\it Feed - back control loop}}
   \vspace{1mm}

\subsection*{3.1 Fractional-order controlled system}

The fractional-order controlled system will be represented with a fractional model with the fractional differential equation given by the 
following expression $(_0D^{\mu}_{t} \equiv D^{\mu}_{t})$:
\begin{equation} \label{n-DR}
     a_{n}\, D^{\beta_{n}}_ty(t) + \ldots +
     a_{1}\, D^{\beta_{1}}_ty(t) +
     a_{0}\, D^{\beta_{0}}_ty(t) =
     b_{n}\, D^{\alpha_{n}}_tu(t) + \ldots +
     b_{1}\, D^{\alpha_{1}}_tu(t) +
     b_{0}\, D^{\alpha_{0}}_tu(t),
\end{equation}
where  $\beta_k, \alpha_k$  $(k = 0, 1, 2, \dots)$ are generally real
numbers, $\beta_{n} > \ldots > \beta_{1} > \beta_{0}$,
$\alpha_{m} > \ldots > \alpha_{1} > \alpha_{0}$ and
$a_k, b_k$ $(k= 0, 1, \dots)$ are arbitrary constants.

For obtaining a discrete model of the fractional-order system
(\ref{n-DR}), we have to use discrete approximations of the
fractional-order integro-differential operators and then we
obtain a general expression for the discrete transfer function of
the controlled system \cite{Vinagre}
\begin{equation}
  G_s(z)= \frac{b_m (\omega(z^{-1}))^{\alpha_m}+ \dots +
          b_1 (\omega(z^{-1}))^{\alpha_1}+b_0(\omega(z^{-1}))^{\alpha_0}}
               {a_n (\omega(z^{-1}))^{\beta_n}+ \dots +
          a_1 (\omega(z^{-1}))^{\beta_1}+a_0(\omega(z^{-1}))^{\beta_0}},
\end{equation}
where $(\omega(z^{-1}))$ denotes the discrete operator, expressed
as a function of the complex variable $z$ or the shift operator
$z^{-1}$.

For discrete time step $k$, according to relation (\ref{DFD})
and the inverse Z-transform of the difference equation (\ref{n-DR}),
we can write the fractional-order difference equation in the form
\begin{eqnarray}\label{r4_a}
  a_n\, T^{-\beta_n}\sum_{j=0}^{k}c_j^{(\beta_n)} y_{k-j} + \dots +
  a_1\, T^{-\beta_1}\sum_{j=0}^{k}c_j^{(\beta_1)} y_{k-j} +
  a_0\, T^{-\beta_0}\sum_{j=0}^{k}c_j^{(\beta_0)} y_{k-j} =
\nonumber \\
  b_m\, T^{-\alpha_m}\sum_{j=0}^{k}c_j^{(\alpha_m)} u_{k-j} + \dots +
  b_1\, T^{-\alpha_1}\sum_{j=0}^{k}c_j^{(\alpha_1)} u_{k-j} +
  b_0\, T^{-\alpha_0}\sum_{j=0}^{k}c_j^{(\alpha_0)} u_{k-j}.
\end{eqnarray}
From relation (\ref{r4_a}) the following form of
the difference equation can be obtained
\begin{equation}
y_k=\frac{\sum_{i=1}^{m}(b_i T^{-\alpha_i}\sum_{j=0}^{k}
          c_j^{(\alpha_i)}u_{k-j})
          - \sum_{i=1}^{n}(a_i T^{-\beta_i}\sum_{j=1}^{k}
          c_j^{(\beta_i)}y_{k-j})}
         {\sum_{i=0}^{n}a_i T^{-\beta_i}c_0^{(\beta_i)}},
\end{equation}
for $k=2, 3, \dots$, where $y_0=0$ and $y_1=0$.

\subsection*{3.2 Fractional order controller}

The fractional $PI^{\lambda}D^{\delta}$ controller will be represented
by discrete transfer function given in the following expression:
\begin{equation}\label{Gr}
      G_{c}(z)=\frac{U(z)}{E(z)}=
      K + \frac{T_i}{(\omega(z^{-1}))^{\lambda}}+T_d(\omega(z^{-1}))^{\delta},
\end{equation}
where $\lambda$ and $\delta$ are arbitrary real numbers
$(\lambda, \delta \geq 0)$, $K$ is the proportional constant,
$T_i$ is the integration constant and $T_d$ is the derivative
constant.

Taking $\lambda=1$ and $\delta=1$, we obtain a classical $PID$
controller. If $\lambda = 0$ and/or $T_i=0$, we obtain
a~$PD^{\delta}$ controller, etc.
All these types of controllers are particular cases of the
$PI^{\lambda}D^{\delta}$ controller, which is more flexible and
gives an opportunity to better adjust the dynamical properties
of the fractional-order control system.

Similar to methods for controlled system, they can be used for
obtaining the difference equation of the fractional-order controller
(\ref{Gr}).

\section*{4. ILLUSTRATIVE EXAMPLE}

We give in this section an example of modelling the
stable dynamical system by using fractional calculus in
discrete domain. The fractional-order control system consists of
the real controlled system with the coefficients:
\begin{equation}\label{sys}
a_2 = 0.8, \,\, a_1 = 0.5, \,\, a_0 = 1.0, \,\,
\beta_2 = 2.2, \,\, \beta_1 = 0.9, \,\, \beta_0 = 0, \,\,
b_0 = 1.0, \,\, \alpha_0 = 0
\end{equation}
and the fractional-order $PD^{\delta}$ controller, designed on the
stability measure $S_t=2.0$ and damping measure $\xi=0.4$, with
the coefficients:
\begin{equation}\label{cont}
K = 50.0, \,\, T_d = 5.326, \,\, \delta = 1.286.
\end{equation}
The fractional-order differential equation of the closed control loop
has the form $(_0D^{\mu}_{t} \equiv D^{\mu}_{t})$:
\begin{equation} \label{r_44}
  {a_2\, D^{\beta_2}_t y(t) + a_1\, D^{\beta_1}_t y(t) + T_d\, D^{\delta}_t y(t)
  + (a_0 + K) y(t)} = K~w(t) + T_d\,  D^{\delta}_t w(t).
\end{equation}
The resulting differential equation (\ref{r_44}) can be rewritten to the fractional difference equation of the closed control loop in general form
via inverse Z-transform of the relation (\ref{DFD}):
\begin{equation} \label{numer_PD}
       y_k =
       \frac{K w_k+
       T_d{T}^{-\delta}\sum\limits_{j=0}^{k}{c_j^{(\delta)} w_{k-j}}-
       a_2{T}^{-\beta_2}\sum\limits_{j=1}^{k}{c_j^{(\beta_2)} y_{k-j}}-
       a_1{T}^{-\beta_1}\sum\limits_{j=1}^{k}{c_j^{(\beta_1)} y_{k-j}}-
       T_d{T}^{-\delta}\sum\limits_{j=1}^{k}{c_j^{(\delta)} y_{k-j}}
            }
            {
       a_2{T}^{-\beta_2} c_0^{(\beta_2)} +
       a_1{T}^{-\beta_1} c_0^{(\beta_1)} +
       T_d{T}^{-\delta} c_0^{(\delta)} +
       (a_0 +K)
            },
\end{equation}
for $(k=1,2,\dots)$, where  $y_0 = 0,\,w_0 = 0,\, w_1 = 0$\,
and \, $w_k = 1$, for $(k=2,3,\dots)$. The binomial coefficients
were calculated according to relation (\ref{b_k}).

\section*{5. CONCLUSION}

The above methods make it possible to model and analyse
(simulate) fractional-order control systems in the discrete
domain and also to realise digital fractional-order controllers.
We have shown the survey of approximation methods
and Z-transform method as good
approximations of the fractional-order operator
$_aD^{\alpha}_{t}$.
These methods for discretisation of fractional calculus are
suitable for realisation and implementation of
fractional-order controllers, because this controllers
are more robust than classical one.

Stability investigation of discrete fractional-order control
system and practical simulation will be studied in further work.

\section*{ACKNOWLEDGEMENTS}

This work was partially supported by grant {\it \,Fractional-order Dynamical Systems and Controllers:
Discrete and Frequency-Domain Models and Algorithms} from
Austrian Institute for East and South-East Europe
and partially supported by grant VEGA 1/7098/20 from the
Slovak Agency for Science.

\end{document}